\documentclass[12pt]{article}
\usepackage{mathrsfs}
\usepackage{latexsym,amsmath,amssymb,amsfonts,epsfig,graphicx,cite,psfrag}
\usepackage{eepic,color,colordvi,amscd,amsthm}
\usepackage{pstcol,pstricks,color}
\newtheorem{theorem}{Theorem}[section]

 \numberwithin{equation}{section}
 \numberwithin{figure}{section}

\def\qed{\hfill \rule{4pt}{7pt}}

\def\pf{\noindent {\it Proof.} }

\begin{document}

\begin{center}
{\large\bf Partially $2$-Colored Permutations and

the Boros-Moll Polynomials}
\end{center}

\vskip 2mm \centerline{William Y.C. Chen$^1$, Sabrina X.M. Pang$^2$,
and Ellen X.Y. Qu$^3$}

\begin{center}
$^{1}$Center for Combinatorics, LPMC-TJKLC\\
Nankai University, Tianjin 300071, P.R. China

$^{2}$College of Mathematics and Statistics\\
Hebei University of Economics and Business\\
Shijiazhuang, Hebei 050061, P.R. China

$^{3}$School of Mathematical Sciences\\
Ocean University of China\\
Qingdao, Shandong 266100, P.R. China

\vskip 2mm
  $^1$chen@nankai.edu.cn  $^2$stpangxingmei@heuet.edu.cn
  $^3$ellenqu@163.com
\end{center}

\begin{center}
{\bf Abstract}
\end{center}

We find a combinatorial setting for
the coefficients of the Boros-Moll polynomials $P_m(a)$ in terms of partially 2-colored permutations. Using this model, we
give a combinatorial proof of a recurrence relation on the coefficients
 of $P_m(a)$. This approach enables us to
 give a combinatorial interpretation of the log-concavity of $P_m(a)$ which was conjectured by Moll and confirmed by Kauers and Paule.

\vskip 3mm


\vskip 3mm



\section{Introduction}

The main objective of this paper is to present a combinatorial
approach to  the log-concavity of the Boros-Moll polynomials. The
Boros-Moll polynomials  $P_m(a)$ arise in the evaluation of a
quartic integral, see \cite{Bor-Mol99a,Bor-Mol99b,Bor-Mol99c,
Bor-Mol01,Bor-Mol04,Mol02}.
Boros and Moll have shown that for any $a>-1$ and any nonnegative
integer $m$,
\begin{equation}
\int_{0}^{\infty}\frac{1}{(x^4+2ax^2+1)^{m+1}}dx=\frac{\pi}{2^{m+3/2}(a+1)^{m+1/2}}P_m(a),
\end{equation}
where
\begin{equation}\label{ds}
P_m(a)=\sum_{j, k}{2m+1\choose 2j}{m-j\choose k}{2k+2j\choose
k+j}\frac{(a+1)^j(a-1)^k}{2^{3(k+j)}}.
\end{equation}
Boros and Moll also derived a single sum
 formula for $P_m(a)$:
  \begin{equation}\label{ss}
P_m(a)=2^{-2m}\sum\limits_{k}2^k{2m-2k\choose m-k}{m+k\choose
k}(a+1)^k,
 \end{equation}
  which implies that the coefficients of $P_m(a)$ are positive. More precisely, let $d_i(m)$ be the coefficient of $a^i$ in $P_m(a)$. Then (\ref{ss}) gives
\begin{equation}\label{eq1.4}
d_{i}(m)=2^{-2m}\sum\limits_{k=i}^m 2^k{2m-2k\choose m-k}{m+k\choose
k}{k\choose i}.
\end{equation}
Several  proofs of the formula (\ref{ss}) can be found in the survey of
Amdeberhan and  Moll \cite{Amd-Mol08}.

Further positivity properties of $P_m(a)$ have been studied
recently. Boros and Moll \cite{Bor-Mol99c} have shown that the
sequence $\{d_i(m)\}_{0\leq i \leq m}$ is unimodal for $m\geq 0$.
Moll conjectured that this sequence is log-concave, that is, for
$m\geq2$ and $1\leq i \leq m-1$,
\begin{equation}\label{eq1.1}
d_{i}^2(m)\geq d_{i-1}(m)d_{i+1}(m).
\end{equation}
This conjecture has been confirmed by Kauers and Paule
\cite{Kau-Pau07} based on recurrence relations. Chen and Xia \cite{Chen-Xia08} have proved a
stronger property of $d_i(m)$, called the ratio monotone property,
which implies both the log-concavity and the spiral property. Moll
\cite{Man09,Mol07}  posed a conjecture that is stronger than the
log-concavity of  $P_m(a)$.  This conjecture
 has been proved by Chen and
Xia \cite{Chen-Xia09}. Chen and Gu \cite{Chen-Gu09}
established the reverse ultra log-concavity of the Boros-Moll
polynomials.

It turns out that the polynomials  $P_m(a)$ are closely related to
combinatorial structures.  The $2$-adic valuation of the numbers
$i!m!2^{m+i}d_i(m)$ has been studied by Amdeberhan, Manna and Moll
\cite{Amd-Man-Mol08}, and Sun and Moll \cite{Sun-Mol08}. By using
reluctant functions and an extension of Foata's bijection, Chen,
Pang and Qu \cite{Chen-Pang-Qu08} have found a combinatorial
derivation of the single sum formula
 (\ref{ss}) from the  double sum formula (\ref{ds}).   For the
special case   $a=1$, we are led to a combinatorial argument for
the  identity
\begin{equation*}\label{mi}
\sum\limits_{k=0}^{m}2^{-2k}{2k \choose k}{2m-k \choose
m}=\sum\limits_{k=0}^{m}2^{-2k}{2k \choose k}{2m+1 \choose 2k}.
\end{equation*}
However, this combinatorial approach
does not seem to  apply to
recurrence relations for $d_i(m)$ or the
log-concavity of $P_m(a)$.

In this paper, we shall consider a variation of the  coefficients
$d_i(m)$, that is,
\begin{equation}\label{D-im}
D_i(m)={2m \choose m-i}m!i!(m-i)!2^id_i(m).
\end{equation}
Then the numbers $D_i(m)$ have a combinatorial interpretation in
terms of partially 2-colored permutations.

  Using this combinatorial setting,
we give an explanation of the following recurrence relation of
$d_i(m)$ derived independently by Kauers and Paule \cite{Kau-Pau07}
and Moll \cite{Mol07}:
\begin{align}\label{eqn3}
i(i+1)d_{i+1}(m)=i(2m+1)d_{i}(m)-(m-i+1)(m+i)d_{i-1}(m).
\end{align}

The  reasoning of the above recurrence relation also implies
a simple combinatorial interpretation of the log-concavity of the Boros-Moll polynomials.

\section{A combinatorial setting for $D_i(m)$}\label{section2}

In this section, we shall give a combinatorial interpretation of  $D_i(m)$ by introducing the structure of partially 2-colored permutations.
 Throughout this paper, we shall adopt the notation
$(x)_n$ for rising factorials, that is, $(x)_0=1$ and for $n>0$,
\[ (x)_n=x(x+1)\cdots(x+n-1).\]
From the expression \eqref{eq1.4} for $d_i(m)$, we have
\allowdisplaybreaks {\small\begin{align*}
d_{i}(m)&=2^{-2m}\sum\limits_{k=i}^m
2^k{2m-2k\choose m-k}{m+k\choose k}{k\choose i} \\
&= 2^{-2m}\sum\limits_{j=0}^{m-i}
2^{j+i}{2m-2i-2j\choose m-i-j}{m+i+j\choose i+j}{i+j\choose i} \\
&= 2^{-2m}\sum\limits_{j=0}^{m-i}
2^{j+i}\frac{(2m-2i-2j)!}{(m-i-j)!(m-i-j)!}
\cdot\frac{(m+i+j)!}{(i+j)!m!} \cdot\frac{(i+j)!}{j!i!}
 \\
&=  2^{-2m}\sum\limits_{j=0}^{m-i}
2^{j+i}\frac{2^{2m-2i-2j}(m-i-j-\frac{1}{2})!}{(m-i-j)!}\cdot
\frac{(m+i+j)!}{(i+j)!m!} \cdot\frac{(i+j)!}{j!i!}.
\end{align*}}
It follows that {\small\begin{align*} m!i!(m-i)!2^id_{i}(m) & =
(m-i)!\sum\limits_{j=0}^{m-i} \left(\frac{1}{2}\right)^{j}
\frac{(m-i-j-\frac{1}{2})!}
{(m-i-j)!}\cdot\frac{(m+i+j)!}{j!},\\
& = \sum\limits_{j=0}^{m-i}{m-i\choose j}
\left(\frac{1}{2}\right)^{j} \left(\frac{1}{2}\right)_{m-i-j}
(m+i+j)!,
\end{align*}}
which yields {\small\begin{equation} \label{eq2.1} D_i(m)={2m\choose
m-i}\sum\limits_{j=0}^{m-i}{m-i\choose j}
\left(\frac{1}{2}\right)^{j}
\left(\frac{1}{2}\right)_{m-i-j}(m+i+j)!.
\end{equation}}

We proceed to give a combinatorial interpretation of $D_i(m)$
according to the expression  \eqref{eq2.1}.
It is well known that $(x)_n$ equals the generating function for
permutations on $[n]$ with respect to the number of cycles. Let $\sigma$ be a permutation on $[n]$.
The weight of $\sigma$ is defined as $x^k$, where $k$ is the number
of cycles in $\sigma$. So $(x)_n$ is the weighted count of
permutations on $[n]$.

 Suppose that
$(A,B,C)$ is a composition of $[2m]=\{1, 2, \ldots, 2m\}$, namely,
any $A$, $B$ and $C$ are disjoint and $A\cup B\cup C=[2m]$, where
$A$, $B$ and $C$ are allowed to be empty. A permutation on $[2m]$
associated with a composition $(A, B, C)$ of $[2m]$ is called a
partially $2$-colored permutation on $[2m]$ if it can be written as
$(\pi|\sigma)$, where $\pi$ is a permutation on $A\cup B$ and
$\sigma$ is a permutation on $C$. We assume that the elements in $A$
are white, the elements in $B$ are black and written in boldface,
while the elements in $C$ are uncolored.

Moreover, we need to use two different representations for the permutations
$\pi$ and $\sigma$ in a partially $2$-colored permutation $(\pi|\sigma)$. To be precise, we shall
write $\pi$ in the one-line notation in the form of a sequence.  For
example,
$5,7,8,2,1,6,4,3$ is the one-line representation
of a permutation.
On the other hand, we shall express $\sigma$ in terms of the cycle decomposition.
 For
instance, the permutation in the above example has cycle decomposition $(1,5)(2,7,4)(3,8)(6)$.

Let $\mathcal{D}_{i}(m)$ denote the set of all partially 2-colored
permutations $(\pi|\sigma)$ on $[2m]$ such that the 2-colored
permutation $\pi$ has $m+i$ black elements. For example, consider
the partially 2-colored permutation
\[({\bf 2},{\bf 12},8,{\bf11},{\bf5},{\bf9},
{\bf7},1,{\bf4},{\bf3}|(6,10))\]
 in $\mathcal{D}_{2}(6)$. Then we have
$A=\{1,8\}$, $B=\{2,3,4,5,7,9,11,12\}$,  and $C=\{6,10\}$. From the
definition, we see that for a partially 2-colored permutation
$(\pi|\sigma)$ in $\mathcal{D}_{i}(m)$, we have  $|A\cup C|=m-i$.

We are now ready to give a combinatorial interpretation of $D_i(m)$. With respect to the
weight a partially
 $2$-colored permutation $(\pi|\sigma)$ in $\mathcal{D}_{i}(m)$, we impose
  the following rules:
\begin{enumerate}
\item[(1)] An element in $A$ is given a weight $\frac{1}{2}$;
\item[(2)] A cycle in $\sigma$ is given a weight $\frac{1}{2}$.
\end{enumerate}
The weight $(\pi|\sigma)$ is defined as the
product of the weights of the white elements and the cycles. In light of the above weight assignment, $D_i(m)$ can be
viewed as a weighted count of partially 2-colored permutations. The
weight of a set $S$ means to be the sum of weights of its elements,
and is denoted by $w(S)$.

\begin{theorem}\label{theorem3} For $m\geq 1$,
$D_{i}(m)$ equals the  weight of $\mathcal{D}_{i}(m)$.
\end{theorem}

\pf  Given a composition $(A,B,C)$ of $[2m]$ such that $|B|=m+i$ and
$|A\cup C|=m-i$. Assume that there are $j$ elements in $A$. It
is clear that there are $m-i-j$ elements in $C$. Now, there are
${2m\choose m-i}$ ways to distribute $2m$ elements into  $B$ and
$A\cup C$. Moreover, there are ${m-i \choose j}$ ways to distribute
$m-i$ elements into   $A$ and $C$.

Consider partially 2-colored permutations in
$\mathcal{D}_{i}(m)$ associated with composition $(A,B,C)$ of $[2m]$. Since $|A \cup
B|=m+i+j$, the  sum of weights of permutations on $A\cup B$ equals
$$\left(\frac{1}{2}\right)^j\cdot (m+i+j)!.$$ The weighted sum of
permutations on $C$ equals $\left(\frac{1}{2}\right)_{m-i-j}.$ This
completes the proof.\qed

\section{Combinatorial proof of the recurrence relation}

Using the interpretation of $D_i(m)$ in terms of
partially $2$-colors permutation, we   give a combinatorial
proof for the following recurrence relation of the coefficients
$d_{i}(m)$ of the Boros-Moll polynomials
\begin{align}\label{eqn1}
i(i+1)d_{i+1}(m)=i(2m+1)d_{i}(m)-(m-i+1)(m+i)d_{i-1}(m).
\end{align}
This recurrence was independently derived by Kauers, Paule \cite{Kau-Pau07} and
Moll \cite{Mol07}.

Utilizing   \eqref{D-im}, the recurrence relation
\eqref{eqn1} can be restated as
\begin{align}\label{equation1}
\frac{1}{2}(m+i+1)D_{i+1}(m)+2(m-i+1)D_{i-1}(m)=(2m+1)D_{i}(m).\end{align}

 To give a combinatorial proof of \eqref{equation1}, we need to introduce some notation.   Let $\mathcal{A}_{i}(m)$ (resp.
$\mathcal{B}_{i}(m)$ and $\mathcal{C}_{i}(m)$) denote the set of all
partially 2-colored permutations $(\pi|\sigma)$ in
$\mathcal{D}_{i}(m)$ such that exactly one element in $A$ (resp. $B$
and $C$) is underlined. Obviously, the four sets
$\mathcal{A}_{i}(m)$, $\mathcal{B}_{i}(m)$, $\mathcal{C}_{i}(m)$ and
$\mathcal{D}_i(m)$ are disjoint. For example,
\[({\bf 2},{\bf 12},8,{\bf11},{\bf5},\underline{{\bf9}},
{\bf7},1,{\bf4},{\bf3}|(6,10))\] is an underlined
 partially $2$-colored permutation
 belonging to $\mathcal{B}_2(6)$. By definition and Theorem \ref{theorem3}, we have
 \begin{eqnarray}
    (m+i)  {D}_i(m)   & = & w(\mathcal{B}_i(m)), \label{d2}\\[3pt]
    (m-i)  {D}_i(m)  & = & w(\mathcal{A}_i(m)\cup \mathcal{C}_i(m)) .  \label{d1}
 \end{eqnarray}

\pf  From \eqref{d2} and \eqref{d1}, we know that
 \begin{eqnarray}\label{md1}
  (m+i+1)D_{i+1}(m)  & = & w(\mathcal{B}_{i+1}(m)),\\[3pt]
  \label{md2}
    (m-i+1)D_{i-1}(m)   &= &
 w(\mathcal{A}_{i-1}(m) \cup
 \mathcal{C}_{i-1}(m)).
\end{eqnarray}
On the other hand, we have
\begin{equation}\label{ad}
(2m+1)D_{i}(m)  =   w( \mathcal{A}_{i}(m)\cup
\mathcal{B}_{i}(m)\cup\mathcal{C}_{i}(m)\cup \mathcal{D}_{i}(m)).
\end{equation}

First, we claim that
\begin{equation}\label{b2}
{1\over 2} w(\mathcal{B}_{i+1}(m)) =  w(\mathcal{A}_{i}(m)).
\end{equation}
Given $(\pi|\sigma)\in \mathcal{B}_{i+1}(m)$ with underlying
composition $(A, B, C)$, where $|B|=m+i+1$ and $|A\cup C|=m-i-1$, by
changing the  underlined black element in $\pi$ to an underlined
white
 element, we obtain an underlined partially
 $2$-colored permutation in $\mathcal{A}_i(m)$.
 Clearly, this operation yields a bijection
 between $\mathcal{B}_{i+1}(m)$ and $\mathcal{A}_i(m)$.
 Since  the weight of a white element equals $1/2$, we obtain (\ref{b2}).
Substituting $i$ with $i-1$ in (\ref{b2}), we get
\begin{equation}\label{ab}
w(\mathcal{B}_{i}(m))=  2 w(\mathcal{A}_{i-1}(m) ).
\end{equation} Hence
(\ref{equation1}) simplifies to the following relation
\begin{equation}\label{ac2}
2w( \mathcal{C}_{i-1}(m)) = w(\mathcal{C}_{i}(m)\cup \mathcal{D}_{i}(m)).
\end{equation}

Assume that  $(\pi|\sigma)\in \mathcal{C}_{i-1}(m)$ is a partially
$2$-colored permutation with underlying composition $(A, B, C)$,
that is, $|B|=m+i-1$, $|A\cup C|=m-i+1$, and $\sigma$ is a
permutation with an underlined element. Suppose that $\sigma$ has
cycle decomposition $C_0, C_1, \ldots, C_r$, where $C_0$ contains
the underlined element. Without loss of generality, we may always
write $C_0$ as $(\underline{i_1}i_2\cdots i_{k})$.  Given
$(\pi|\sigma) \in \mathcal{C}_{i-1}(m)$, we define
\[ \Delta(\pi|\sigma)=\{\Delta_1, \Delta_2, \ldots, \Delta_k\},\]
where
\begin{eqnarray*}
\Delta_1&=&(\pi, {\mathbf{i_1}}|(\underline{i_2}, i_3, \ldots,
i_k)C_1C_2\cdots C_r),\\ \Delta_2&=&(\pi, {\mathbf{i_1}},
i_2|(\underline{i_3}, \ldots, i_k)  C_1C_2\cdots
C_r),\\
&&\qquad \qquad \quad \cdots\\ \Delta_{k-1}&=&(\pi, {\mathbf{i_1}},
i_2,
\ldots, i_{k-1}|(\underline{i_k} ) C_1C_2\cdots C_r),\\
\Delta_k&=&(\pi, {\mathbf{i_1}}, i_2, \ldots, i_{k-1}, i_{k}|C_1C_2\cdots
C_r).
\end{eqnarray*}

For $1\leq j\leq k-1$,  we have $\Delta_j\in
\mathcal{C}_{i}(m)$ and
\begin{equation}\label{w1}
w(\Delta_j)= {1\over 2^{j-1}} w(\pi|\sigma).
\end{equation}Moreover, we see that $\Delta_k \in
\mathcal{D}_{i}(m)$ and
\begin{equation}\label{w2}
w(\Delta_k)= {1\over 2^{k-2}} w(\pi|\sigma).
\end{equation}
Conversely, any partially colored permutation in
$\mathcal{C}_{i}(m)\cup  \mathcal{D}_{i}(m)$ can be obtained from a partially colored permutation
in $ \mathcal{C}_{i-1}(m)$ by applying the above
operation $\Delta$. Thus, we deduce that
\begin{equation} \label{cd}
 \Delta (\mathcal{C}_{i-1}(m)) = \mathcal{C}_{i}(m)\cup  \mathcal{D}_{i}(m),
\end{equation}
where $\Delta$ acts on the partially colored permutations in
$\mathcal{C}_{i-1}(m)$. Since
\begin{equation*}\sum_{j=1}^{k-1}\frac{1}
{2^{j-1}}+\frac{1}{2^{k-2}}=2,\end{equation*} combining (\ref{w1}),
(\ref{w2}) and (\ref{cd}) we obtain (\ref{equation1}). This
completes the proof. \qed

\section{ Combinatorial proof of the log-concavity}\label{section4}

In this section, we shall use the structure of
partially $2$-colored permutations to give a combinatorial reasoning of the following relation \begin{eqnarray}\label{eq1.5}
(m+i+1) D_{i+1}(m)\cdot(m-i+1) D_{i-1}(m) <(m+i)(m-i+1)D^2_{i}(m),
\end{eqnarray}
which implies the log-concavity of the Boros-Moll polynomials. We
shall follow the notation introduced in the previous section.

\pf  From (\ref{md1}) and (\ref{md2}), we see that
\begin{eqnarray}
 \lefteqn{(m+i+1) D_{i+1}(m)\cdot(m-i+1)
D_{i-1}(m)}\nonumber \\[3pt]
&=& w(\mathcal{B}_{i+1}(m))\cdot w(\mathcal{A}_{i-1}(m) \cup
 \mathcal{C}_{i-1}(m))\nonumber\\[3pt]
 &=& w(\mathcal{B}_{i+1}(m))\cdot w(\mathcal{A}_{i-1}(m))+w(\mathcal{B}_{i+1}(m))\cdot
 w(\mathcal{C}_{i-1}(m)).\label{e1}
  \end{eqnarray}
Meanwhile, in view of \eqref{d2} and \eqref{d1}, we find
\begin{eqnarray}
\lefteqn{(m+i)(m-i+1)D^2_{i}(m)}\nonumber \\[3pt]
& = & w(\mathcal{B}_{i}(m))\cdot
w(\mathcal{A}_{i}(m) \cup
 \mathcal{C}_{i}(m)
 \cup\mathcal{D}_i(m))\nonumber\\[3pt]
 &=& w(\mathcal{B}_{i}(m))\cdot w(\mathcal{A}_{i}(m))+
 w(\mathcal{B}_{i}(m))\cdot w(\mathcal{C}_{i}(m)\cup \mathcal{D}_{i}(m)).  \label{e2}
 \end{eqnarray}
Hence (\ref{eq1.5}) can be recast as
\begin{eqnarray}
&&w(\mathcal{B}_{i+1}(m))\cdot
w(\mathcal{A}_{i-1}(m))+w(\mathcal{B}_{i+1}(m))\cdot
 w(\mathcal{C}_{i-1}(m))\nonumber\\[3pt]
 &&\qquad\qquad<w(\mathcal{B}_{i}(m))\cdot w(\mathcal{A}_{i}(m))+
 w(\mathcal{B}_{i}(m))\cdot w(\mathcal{C}_{i}(m)\cup
 \mathcal{D}_{i}(m)).\label{e3}
\end{eqnarray}
Invoking \eqref{b2} and \eqref{ab},  we obtain
\begin{eqnarray} \label{ba}
w(\mathcal{B}_{i+1}(m))\cdot
w(\mathcal{A}_{i-1}(m))&=&w(\mathcal{B}_{i}(m))\cdot
w(\mathcal{A}_{i}(m)).
\end{eqnarray}
Using (\ref{ba}) and the fact that
 \[ 2w( \mathcal{C}_{i-1}(m))  =  w(\mathcal{C}_{i}(m)\cup
\mathcal{D}_{i}(m))\]  as given by (\ref{ac2}), \eqref{e3}
simplifies to
\begin{eqnarray}\label{e4}
{1\over 2}w(\mathcal{B}_{i+1}(m))<w(\mathcal{B}_{i}(m)).
\end{eqnarray}
Applying(\ref{b2}), (\ref{e4}) is
equivalent to the relation
\begin{equation} \label{a-b}
w(\mathcal{A}_{i}(m))<w(\mathcal{B}_{i}(m)),
\end{equation}
which can be easily deduced from \eqref{d2} and \eqref{d1},
since for $1\leq i\leq m-1$,
\begin{eqnarray}
w(\mathcal{A}_{i}(m))\leq w(\mathcal{A}_{i}(m) \cup \mathcal{C}_i(m)) = (m-i) {D}_i(m) <(m+i) {D}_i(m) =w(\mathcal{B}_{i}(m)).
\end{eqnarray}
This completes the proof. \qed

\vskip 3mm \noindent {\bf Acknowledgments.} This work was supported
by the 973 Project, the National Natural Science Foundation of
China, the PCSIRT Project  and the Fundamental Research Funds for
Central Universities of the Ministry of Education of China.


\begin{thebibliography}{99}

\bibitem{Amd-Man-Mol08} T. Amdeberhan, D. Manna and V. Moll, The
$2$-adic valuation of a sequence arising from a rational integral,
J. Combin. Theory Ser. A 115 (8) (2008) 1474--1486.

\bibitem{Amd-Mol08} T. Amdeberhan and V. Moll, A formula for a quartic integral:
A survey of old proofs and some new ones, Ramanujan J. 18 (2009)
91--102.

\bibitem{Bor-Mol99a} G. Boros and V. Moll, An integral hidden in
Gradshteyn and Ryzhik, J. Comput. Appl. Math. 106 (1999) 361--368.

\bibitem{Bor-Mol99b} G. Boros and V. Moll,
A sequence of unimodal polynomials, J. Math. Anal. Appl. 237 (1999)
272--287.

\bibitem{Bor-Mol99c}G. Boros and V. Moll, A criterion for unimodality,
Electron. J. Combin. 6 (1999) \#R10.

\bibitem{Bor-Mol01} G. Boros and V. Moll,
The double square root, Jacobi polynomials and Ramanujan's Master
Theorem, J. Comput. Appl. Math. 130 (2001) 337--344.

\bibitem{Bor-Mol04}
G. Boros and V. Moll, Irresistible Integrals, Cambridge University
Press, New York/Cambridge, 2004.

\bibitem{Chen-Gu09} W.Y.C. Chen and C.C.Y Gu,
The reverse ultra log-concavity of the Boros-Moll polynomials, Proc.
Amer. Math. Soc. 137 (2009) 3991--3998.



\bibitem{Chen-Pang-Qu08}
W.Y.C. Chen, S.X.M. Pang and E.X.Y. Qu, On the combinatorics of the
Boros-Moll polynomials, Ramanujan J. 21 (2010) 41--51.

\bibitem{Chen-Xia08} W.Y.C. Chen and E.X.W. Xia,
The ratio monotonicity of the Boros-Moll polynomials, Math. Comp.
78 (2009) 2269--2282.

\bibitem{Chen-Xia09} W.Y.C. Chen and E.X.W. Xia, A proof of Moll's minimum conjecture,
European J. Combin., to appear.



\bibitem{Kau-Pau07} M. Kauers and P. Paule,
A computer proof of Moll's log-concavity conjecture, Proc. Amer.
Math. Soc. 135 (2007) 3847--3856.


\bibitem{Mol02} V. Moll, The evaluation of integrals: A personal
story, Notices Amer. Math. Soc. 49 (3) (2002) 311--317.

\bibitem{Mol07} V. Moll,
Combinatorial sequences arising from a rational integral, Online J.
Anal. Comb. 2 (2007) \#4 .

\bibitem{Man09}
V.H. Moll and D.V. Manna, A remarkable sequence of integers, Expo.
Math. 27 (2009) 289--312.

\bibitem{Sun-Mol08} X.Y. Sun and V. Moll,
A binary tree representation for the $2$-adic valuation of a
sequence arising from a rational integral, Integers 10 (2009)
211--222.

\end{thebibliography}
\end{document}